\documentclass[dvipsnames]{amsart}
\usepackage{amsthm}
\usepackage{amssymb}
\usepackage{amsmath}
\usepackage{mathtools} 

\usepackage[english]{babel}
\usepackage{tikz}
\usepackage[shortlabels]{enumitem}
\usepackage{csquotes}

\usepackage[style=alphabetic, backend=biber, sorting=nyt, doi=false, isbn=false,url=false, hyperref,backref,backrefstyle=none, maxnames=20, maxalphanames=20]{biblatex}
\newbibmacro{string+doiurlisbn}[1]{%
  \iffieldundef{doi}{%
    \iffieldundef{url}{%
      \iffieldundef{isbn}{%
        \iffieldundef{issn}{%
          #1%
        }{%
          \href{https://books.google.com/books?vid=ISSN\thefield{issn}}{#1}%
        }%
      }{%
        \href{https://books.google.com/books?vid=ISBN\thefield{isbn}}{#1}%
      }%
    }{%
      \href{\thefield{url}}{#1}%
    }%
  }{%
    \href{https://dx.doi.org/\thefield{doi}}{#1}%
  }%
}
\DeclareFieldFormat{title}{\usebibmacro{string+doiurlisbn}{\mkbibemph{#1}}}
\DeclareFieldFormat[article,incollection,inproceedings,unpublished,misc,book]{title}%
    {\usebibmacro{string+doiurlisbn}{\mkbibquote{#1}}}
\DefineBibliographyStrings{english}{%
  backrefpage = {$\uparrow$},%
  backrefpages = {$\uparrow$}%
}
\usepackage{url}
\usepackage[hyperindex,breaklinks]{hyperref}

\theoremstyle{plain}

\theoremstyle{definition}

\numberwithin{equation}{section}

\def\Q{\mathbb{Q}}
\def\Z{\mathbb{Z}}
\def\F{\mathbb{F}}
\def\R{\mathbb{R}}

\DeclareMathOperator{\Log}{Log}
\DeclareMathOperator{\ICM}{ICM}
\DeclareMathOperator{\Pic}{Pic}

\DeclareMathOperator{\Tr}{Tr}

\newcommand{\p}{{\mathfrak p}}
\newcommand{\q}{{\mathfrak q}}
\newcommand{\frP}{{\mathfrak P}}
\newcommand{\cB}{{\mathcal B}}
\newcommand{\cO}{{\mathcal O}}

\newcommand{\cS}{{\mathcal S}}
\renewcommand{\bar}{\overline}
\newcommand{\vphi}{{\varphi}}
\newcommand{\set}[1]{\left\lbrace#1\right\rbrace }
\newcommand{\type}{\mathrm{type}} 

\usepackage{color}
\newcommand{\smm}[1]{\marginpar{\parbox{3.5cm}{\scriptsize \color{blue} \sf SM:\  #1}}} 

\title{Labeling abelian varieties over finite fields}
\date{\today}

\author[E. Costa]{Edgar Costa}
\address{Edgar Costa, Department of Mathematics, Massachusetts Institute of Technology, Cambridge, MA 02139-4307, USA}
\email{edgarc@mit.edu}
\urladdr{\url{https://edgarcosta.org}}

\author[T. Dupuy]{Taylor Dupuy}
\address{Taylor Dupuy, Department of Mathematics and Statistics,
	Innovation Hall E220,
	82 University Place,
	Burlington, VT, 05405, USA}
\email{taylor.dupuy@gmail.com}
\urladdr{\url{http://tdupu.github.io}}

\author[S. Marseglia]{Stefano Marseglia}
\address{Stefano Marseglia, Laboratoire Jean Alexandre Dieudonné, Université Côte Azur, 06108 Nice Cedex 2, France}
\email{stefano.marseglia@univ-cotedazur.fr}
\urladdr{\url{https://stmar89.github.io}}

\author[D. Roe]{David Roe}
\address{David Roe, Department of Mathematics, Massachusetts Institute of Technology, Cambridge, MA 02139-4307, USA}
\email{roed@mit.edu}
\urladdr{\url{https://math.mit.edu/~roed}}

\author[C. Vincent]{Christelle Vincent}
\address{Christelle Vincent, Department of Mathematics and Statistics,
	Innovation Hall E220,
	82 University Place,
	Burlington, VT, 05405, USA}
\email{christelle.vincent@uvm.edu}
\urladdr{\url{https://www.uvm.edu/~cvincen1}}

\begin{document}

\begin{abstract}
    We describe a deterministic process to associate a practical, permanent label to isomorphism classes of abelian varieties defined over finite fields with commutative endomorphism algebra as long as they are ordinary or defined over a prime field. In the ordinary case, we also provide labels for the polarizations they admit.
\end{abstract}

\maketitle

\section{Introduction}
Let $k$ be a field of size $q=p^a$, where $p$ is a 
prime.
Building on the work of Deligne \cite{Del69}, Howe \cite{Howe95} and Centelghe--Stix \cite{CentelegheStix15} 
the third author developed in \cite{MarsegliaICM, MarsegliaAbVar, MarsegliaCohenMacaulay,Marseglia_local_isom} tools that, when combined with the work of Hofmann--Sircana \cite{HofmannSircana}, allow the enumeration of $k$-isomorphism classes of abelian varieties of dimension $g$ defined over $k$ contained in isogeny classes satisfying \eqref{assumption1} and \eqref{assumption2}, as well as their polarizations when the class is ordinary.
Along with Mckenzie West, the authors have implemented these algorithms and enumerated some of these isomorphism classes and low-degree polarizations
for forthcoming inclusion in the LMFDB \cite{LMFDB} to expand the data about isogeny classes (see \cite{LMFDB_isogeny_paper}).

The purpose of this note is to describe a deterministic process to associate a \textbf{label}, by which we mean a unique permanent identifier, to each of these objects. We believe that useful, short labels that can gain widespread adoption in the community are crucial for the long-term usefulness and intelligibility of results in the field. 
Note that, unlike the labeling scheme for isogeny classes, the scheme we propose does not include enough data to directly recompute the abelian variety (without re-enumerating the isogeny class); it is also not stable under base change.

The third author's algorithms rely on a bijection between certain isomorphism classes of ideals and the isomorphism classes of abelian varieties to be enumerated (\S \ref{sec:theory}). Accordingly, the labeling scheme we propose ultimately labels these isomorphism classes of ideals. Furthermore, when the abelian variety is ordinary we also compute polarizations of low degree; this is done by defining a distinguished representative $I$ in the isomorphism class of ideals and giving a polarization $\lambda$ on the corresponding variety $A$ as an element of the endomorphism algebra.

We close by describing our proposed labeling: a polarized abelian variety corresponding to a pair $(I,\lambda)$ as above is labeled as 
\begin{equation}\label{E:labeling-scheme1}
\texttt{g.q.isog\! \textsf{-}\! N.i.w.j\! \textsf{-}\! d.k},
\end{equation}
where \texttt{g.q.isog} is the label of the isogeny class (\texttt{g} is the dimension, \texttt{q} is the size of the finite field, \texttt{isog} encodes the coefficients of the Weil polynomial; see \cite{lmfdb_permalink_isog_label}), 
\texttt{N.i.w.j} is the label of the isomorphism class of $I$ (\S\S \ref{subsec:isom}, \ref{subsec:weakequivclasses} and \ref{subsec:Pic} -- \texttt{N.i} determines the endomorphism ring \S \ref{subsec:overorders}),
\texttt{d} is the degree $d$ of $\lambda$ and \texttt{k} is determined by the sort key of $\lambda$ (\S \ref{subsec:polarizations}). When no polarization $\lambda$ is given the last part is omitted.

\section{Abelian varieties and ideal classes}\label{sec:theory}

Throughout, let $\mathcal{I}$ be a $k$-isogeny class of abelian varieties of dimension $g$ 
satisfying:
\begin{align}
&\text{the $k$-endomorphism ring of any abelian variety in $\mathcal{I}$ is commutative and} \label{assumption1}\\
&\text{the isogeny class is ordinary or $k$ is the prime field $\F_p$.} \label{assumption2}
\end{align}
By \cite{Del69} and \cite{CentelegheStix15}, there is an equivalence between the category of abelian varieties in $\mathcal{I}$ (with $k$-homomorphisms) and the category of fractional $R$-ideals (with $R$-linear morphisms), where $R$ is the \emph{Frobenius order}, a ring attached to $\mathcal{I}$. 
This equivalence induces a bijection between $\ICM(R)$, the \emph{ideal class monoid} of $R$, and the isomorphism classes in $\mathcal{I}$, see \cite{MarsegliaAbVar}. In this section we introduce the related notation, definitions, and results we need, and refer to \cite{MarsegliaCohenMacaulay,Marseglia_local_isom} for details.

Let $h(x)$ be the characteristic polynomial of any abelian variety in $\mathcal{I}$. By \cite[Theorem 2.(c)]{Tate66} and assumpion \eqref{assumption1}, the $2g$ complex roots of $h(x)$ are distinct.
We sort the distinct monic irreducible factors $h_1(x),\ldots,h_n(x)$ over $\Q$ of $h(x)$ according to the lexicographical order of their coefficients, starting from the constant term.

Consider the étale $\Q$-algebra $K\coloneqq \Q[x]/h(x)$;
throughout we write
\begin{equation}\label{eq:K_as_prod}
    K = K_1 \times \cdots \times K_n,
\end{equation}
where $K_i\coloneqq \Q[x]/h_i(x)$ is a number field, and denote by $F$ the class of the variable $x$ in $K$ and by $\cB_K$ the ordered $\Q$-basis
$(V^{g-1},\ldots,V,1,F,\ldots,F^g)$ of $K$,
where $V=q/F$.
Recall that an \textbf{order} $S$ in $K$ is a subring of $K$ which is a full $\Z$-lattice, that is, the underlying additive group of $S$ has rank $2g$. 
The \textbf{Frobenius order} is then $R \coloneqq \Z[F,V]$, the $\Z$-span of $\cB_K$, 
and an \textbf{overorder} (of $R$) is an order in $K$ containing $R$. We denote the unique maximal order of $K$ by $\cO_K$; using \eqref{eq:K_as_prod} we have
$ \cO_K = \cO_{K_1} \oplus \cdots \oplus \cO_{K_n},$
where $\cO_{K_i}$ is the ring of integers of $K_i$.

For $S$ an order in $K$, a \textbf{fractional $S$-ideal} is a full $\Z$-lattice in $K$ which is closed under multiplication by elements of $S$. A fractional $S$-ideal $I$ is \textbf{invertible} if for every maximal ideal $\p$ of $S$ there exists $a \in K^{\times}$ such that $I_\p = a S_\p$. If $S$ is not maximal then not all fractional ideals are invertible, and ideal multiplication induces a commutative monoid structure on the set of all fractional ideals of $S$.

The \textbf{multiplicator ring} of a fractional $R$-ideal $I$ is the overorder defined by $(I:I)\coloneqq \set{ a\in K: aI\subseteq I}$.
Two fractional $R$-ideals $I$ and $J$ are \textbf{weakly equivalent} if $I_\p\simeq J_\p$ as $R_\p$-modules for each maximal ideal $\p$ of $R$.
The multiplicator ring is an invariant of the weak equivalence class, and if $S$ is an overorder, we denote by $W_S$ the set of weak equivalence classes of fractional $R$-ideals with multiplicator ring $S$.

We need a finer equivalence relation on fractional $R$-ideals: two ideals $I$ and $J$ are \textbf{isomorphic} if there exists $a\in K^\times$ such that $I=aJ$.
Multiplication is well defined on these classes, and the set of isomorphism classes $[I]$ of fractional $R$-ideals $I$ is denoted $\ICM(R)$, the \textbf{ideal class monoid of $R$}. 
 Within $\ICM(R)$, the set of isomorphism classes of invertible fractional $S$-ideals for $S$ an overorder is denoted by $\Pic(S)$; ideal multiplication induces a group structure on $\Pic(S)$.

The \textbf{trace dual ideal} of a fractional $R$-ideal $I$ is $I^t \coloneqq \set{  a \in K : \Tr_{K/\Q}(a I) \subseteq \Z }$. 
For $S$ an overorder and $\mathfrak{p}$ a maximal ideal of $S$, we denote by $\type_{\p}(S)$ the \textbf{Cohen-Macaulay type of $S$ at $\p$}, defined as $\dim_{(S/\p)}(S^t/\p S^t)$.
Note that $\type_{\p}(S)=1$ whenever $\p$ is invertible.
We then define the \textbf{Cohen-Macaulay type} of $S$ to be $\type(S)\coloneqq  \max_{\p}\{\type_{\p}(S)\}$ as $\p$ ranges over the maximal ideals of $S$.

When $\mathcal{I}$ is ordinary, we can compute the polarizations on an abelian variety $A$ belonging to $\mathcal{I}$ from the data of the corresponding fractional $R$-ideal $I$ in the following manner: 
let $\Phi \coloneqq \set{\vphi_1,\ldots,\vphi_g}$ be a CM-type of $K$ which satisfies the Shimura-Taniyama condition; see for example \cite[$\S$ 2.1.4.1]{CCO}.
By \cite{Howe95}, the dual abelian variety $A^\vee$ corresponds to $\bar{I}^t$, where $\bar{\cdot}$ denotes 
complex conjugation on $K$, 
defined by $\bar{F}\coloneq V$, 
and a polarization 
corresponds to an element $\lambda \in K^\times$ such that $\lambda I\subseteq \bar{I}^t$, $\lambda=-\bar\lambda$ and $\Im(\vphi(\lambda))>0$ for every $\vphi \in \Phi$.
Moreover, its degree is $\vert \bar{I}^t/\lambda I\vert$.
Two polarizations $\lambda$ and $\lambda'$ of $I$ are isomorphic if there exists $u\in S^\times$ such that $\lambda = u\bar u \lambda'$.

\section{Labeling ideal classes and polarizations}\label{sec:labeling}
The goal of this section is to define a deterministic procedure to attach a label and a distinguished representative to each ideal class in $\ICM(R)$, where $R$ is the Frobenius order defined in $\S 2$, and to the polarizations they admit, when applicable.

\subsection{Labeling overorders}\label{subsec:overorders}
An overorder $S$ is labeled $\texttt{N.i}$ where $N=[\cO_K:S]$ and $i$ is the index of $S$ when all overorders of index $N$ in $\cO_K$ are sorted in lexicographic order according to the sort key $s(S)\coloneqq [d,n_1,\ldots,n_{g(2g+1)}]$, defined as follows. 
Let $H$ be the matrix whose rows are the coefficients of any $\Z$-basis of $S$ written with respect to $\cB_K$. 
Then $d$ is the least common multiple of the denominators of the entries of $H$ and $n_1,...,n_{g(2g+1)}$ are the entries of the upper triangular part of the (upper triangular) Hermite Normal Form of $dH$.

\subsection{Sorting maximal ideals}\label{subsec:max_idls}
In what follows we will need ordered sets of maximal ideals for a fixed overorder $S$.
We do this by defining sort keys and using the lexicographical order; while our process only sorts maximal ideals within the same order, this is sufficient for our purposes.
Let $p_i \colon K\to K_i$ be the natural projection onto the $i$-th component in \eqref{eq:K_as_prod}.
For $\frP$ a maximal ideal of $\cO_K$, there exists a unique index $1\leq j \leq n$ such that $p_j(\frP)$ is a maximal ideal of $\cO_{K_j}$ and $p_l(\frP) = \cO_{K_l}$ for $l\neq j$.
We define then the sort key $s(\frP)$ of $\frP$ to be $[j,m,n]$ where $\texttt{m.n}$ is the LMFDB label of $p_j(\frP)$; here $m$ is the norm of the ideal and $n$ is a tiebreaker \cite{CremonaPageSutherland}.
If $\p$ is a maximal ideal of a non-maximal order $S$, then the sort key 
of $\p$ is defined to be the lexicographically smallest sort key among those of the finitely many maximal ideals $\frP$ of $\cO_K$ above $\p$.

\subsection{Labeling and distinguished representatives of ideal classes}\label{subsec:isom}
Let $J$ be a fixed fractional $R$-ideal with multiplicator ring $S$.
For any fractional $R$-ideal $I$ which is weakly equivalent to $J$ we have $I=J(I:J)$ where $(I:J)\coloneqq \set{a\in K : aJ \subseteq I}$ is an invertible fractional $S$-ideal.
In fact, by \cite[Corollary~4.5, Theorem~4.6]{MarsegliaICM}, we have a bijection
\[ \ICM(R) \longleftrightarrow \bigsqcup_{R\subseteq S \subseteq \cO_K} \left( W_S \times \Pic(S) \right),\qquad [I]\mapsto (\omega_J,[(I:J)]) \]
where $S$ runs over the finitely many overorders of $R$.

We define the distinguished representative of an arbitrary ideal class $[I]$ as 
the multiplication of the distinguished representative $J$ of the weak equivalence class $\omega_I$ of $I$ (see \S \ref{subsec:weakequivclasses}) with the distinguished representative of $[(I:J)]$ (see \S \ref{subsec:Pic}).
We define the label \texttt{N.i.w.j} of an arbitrary ideal class $[I]$ as the concatenation of the label \texttt{N.i.w} of $\omega_I$ defined in \S \ref{subsec:weakequivclasses} together with the index \texttt{j} of $[(I:J)]\in \Pic(S)$ when enumerated using the sort-key defined in \S \ref{subsec:Pic}.

\subsection{Labeling and distinguished representatives of weak equivalence classes}\label{subsec:weakequivclasses}
We now construct a label and a distinguished representative $J$ for each weak equivalence class $\omega$ of $W_S$ for a fixed overorder $S$.
Throughout, let $\texttt{N.i}$ be the label of $S$ and $\cS \coloneqq (\p_1,\ldots,\p_u)$ be the ordered set of non-invertible maximal ideals of $S$.
Define the label (resp.~distinguished representative) of the invertible class of $W_S$ as \texttt{N.i.1} (resp.~$1\cdot S$).
For each class $\omega \in W_S$ and any $I \in \omega$, set $s(\omega) \coloneqq [1]$ if $\cS = \emptyset$ and 
\[ s(\omega) \coloneqq \left[ \dim_{(S/\p_i)} (I/\p_i I) : 1 \leq i \leq u \right], \]
otherwise. 
If $\type(S)\leq 2$, the string $s(\omega)$ is a complete invariant: 
if $\type_{\p_i}(S)=1$ then $I_{\p_i} \simeq S_{\p_i}$ by, for example, \cite[Proposition~3.4]{MarsegliaCohenMacaulay}; if $\type_{\p_i}(S)=2$ then either $I_{\p_i} \simeq S_{\p_i}$ or $I_{\p_i} \simeq S_{\p_i}^t$ by \cite[Theorem~6.2]{MarsegliaCohenMacaulay}.
We lexicographically sort $W_S$ accordingly and define the label of the $w$th class as $\texttt{N.i.w}$.

If $\type(S)=1$, then we are done, since $W_S$ consists only of the invertible class.
If $\type(S)=2$, it remains to define the distinguished representative of each non-invertible class $\omega\in W_S$:
Let $d$ be the smallest positive integer such that $d S^t \subseteq S$, $\cS_0\coloneqq (\q_1,\ldots,\q_r) \subseteq \cS$ be the ordered set of maximal ideals of $S$ at which $\type_{\q_i}(S)=2$, and $m_1,\ldots,m_r$ be positive integers such that $\q_i^{m_i}S_{\q_i} \subseteq (dS^t)_{\q_i}$.
Then for $I$ any member of $\omega$, we define the distinguished representative $J$ of $\omega$ to be
\[ J \coloneqq \sum_{i=1}^r ( (I_i + \q_i^{m_i} )\prod_{j \ne i}\q_j^{m_j} ), \]
where $I_i \coloneqq S$ if $I_{\q_i}\simeq S_{\q_i}$ and $I_i \coloneqq dS^t$ if $I_{\q_i}\simeq (dS^t)_{\q_i}$; see \cite[Lemma~6.4]{MarsegliaCohenMacaulay}.

Now, we consider the case $\type(S) > 2$.
The following procedure is more time consuming than the previous one.
For each $\p_i \in \cS$, put $T_i\coloneqq (\p_i:\p_i)$; $S \subsetneq T_i$ yields a surjective group homomorphism $\Pic(S) \to \Pic(T_i)$ induced by the extension map.
Let $\mathfrak{K}_i$ be its kernel and set $G_i \coloneqq (T_i^\times/S^\times) \times \mathfrak{K}_i$.
Sort the orders $T_i$ by the size of $G_i$ from smallest to largest, breaking ties using the ordering on the ideals $\p_i$. 
Let $T$ be first among the sorted orders $T_i$.
Let $U$ be a transversal of $T^\times/S^\times$ and $\mathcal K$ be a set of representatives $L$ of the corresponding $\mathfrak{K}_i$ satisfying $LT=T$. By recursion, we assume that we have already computed distinguished representatives for the elements of $W_T$, which we denote by $J_1,...,J_t$.

Fix now a non-invertible $\omega\in W_S$.
By \cite[Proposition~6.2]{Marseglia_local_isom}, the class $\omega$ admits a representative $I_0$ such that $I_0T=J_i$ for a unique index $i$.
Every fractional $S$-ideal $I_1$ weakly equivalent to $I_0$ satisfying $I_1T=J_i$ is of the form $I_1=u\cdot L \cdot I_0$ for unique $u \in U$ and $L \in \mathcal K$.
Since $U$ and $\mathcal{K}$ are finite sets, we list all such ideals $I_1$ and sort them according to their sort key $s(I_1)$, defined in the same way as if $I_1$ were an order $(\S 3.1)$.
Finally, define the distinguished representative $J$ of $\omega$ to be the first ideal of the sorted list
and the sort key $\tilde{s}(\omega)$ of $\omega$ as $s(\omega)$ concatenated with $s(J)$.
Sort $W_S$ accordingly and define the label of the $w$th class as $\texttt{N.i.w}$.

\subsection{Labeling and distinguished representatives of invertible ideal classes}\label{subsec:Pic}\
We start by computing some description of $\Pic(R)$ using \cite{KP}; we then need to fix an ordering of the elements and ideals of $R$ representing each.  We choose a set $\mathcal{P}$ of generators of $\Pic(R)$ by iterating over the maximal ideals of $R$, sorted first by norm and then by the sort-key defined in \S \ref{subsec:max_idls}, keeping only the ones that enlarge the group generated, until we generate all of $\Pic(R)$.  
Now, consider an overorder $S$, and
write $\Pic(S) \simeq \Z/m_1\Z\times\cdots\times\Z/m_k\Z$ with $m_i \mid m_{i+1}$.
The ordered set $\mathcal{P}_S\coloneqq ( \p S : \p \in \mathcal{P} )$ generates $\Pic(S)$, and we seek to construct a basis from it by iteratively choosing elements $g_k, \dots, g_1 \in \Pic(S)$ of order $m_k, \dots, m_1$. 
At stage $i$, let $H_i \coloneqq \langle g_{i+1}, \dots, g_k \rangle$ and choose $L_i$ with class $g_i \in \Pic(S)$ to be the first product of elements of $\mathcal{P}_S$ (by lexicographic order on exponent vectors) with order $m_i$ both in $\Pic(S)$ and $\Pic(S) / H_i$. 
We sort the elements $g\in\Pic(S)$ by writing them as $g=g_1^{e_1}\cdots g_k^{e_k}$ with $0\leq e_i \leq m_i-1$ and we define the distinguished representative of $g$ as $L_1^{e_1}\cdots L_k^{e_k}$.

\subsection{Labeling and distinguished representative of polarizations}\label{subsec:polarizations} 
Assume now that the isogeny class $\mathcal{I}$ is ordinary, and that $I$ is the distinguished representative of its isomorphism class.
Set $S\coloneqq (I:I)$.
For each isomorphism class of polarizations of $I$, we define a distinguished representative as follows.
We will sort elements o$a\in K^\times$ by lexicographically ordering the sequence $(e, n_1, \dots, n_{2g})$ where $(\frac{n_1}{e}, \dots, \frac{n_{2g}}{e})$ are the coefficients of $a$ written with respect to the basis $\cB_K$.
The image of the multiplicative group $U \coloneqq \langle u\bar u\ :\ u \in S^\times \rangle$ under the map
\[ \Log_{\Phi}:K \longrightarrow \R^g \quad  a \mapsto (\log\vert\vphi_1(a)\vert,\ldots,\log\vert\vphi_g(a)\vert),\]
is a lattice $\mathcal{L}$ in $\R^g$.
Fix a representative $\lambda$ of a fixed isomorphism class of polarizations of $I$. 
Consider the set $L$ of elements $u \in S^\times$ giving the points $\Log_{\Phi}(u\bar u) \in \mathcal{L}$ which are the closest to $\Log_{\Phi}(\lambda)\in \R^g$.
For each $u\in L$, add $uv$ to $L$ where $v$ runs over the units in $S^\times$ such that $\vert \Log_{\Phi}(u\bar u\lambda)\vert = \vert\Log_{\Phi}(uv\bar{uv}\lambda)\vert$,
so that $L$ is independent of the choice of the initial representative $\lambda$.
Let $u_0$ be the first element of $L$ according to sorting scheme described above.
Then we set the distinguished representative of the isomorphism class of $\lambda$ to be $\lambda_0 \coloneqq \lambda u_0\bar u_0$, and sort the isomorphism classes of polarizations of $I$ of the same degree $d$ by sorting their distinguished representatives using the procedure described above.

\subsection{Examples}
A Magma implementation of the procedures described above is available at \cite{AbVarFq_LMFDBLabels}.
We list some interesting examples.
\begin{enumerate}[(i)]
    \item 
    For $g=2$, the smallest $q$ so that there exists an isogeny class with an endomorphism ring with Cohen Macaulay type $3$ (largest possible for $g=2$ by \cite[Proposition~4.9]{MarsegliaCohenMacaulay}) has label 
    \href{https://abvar.lmfdb.xyz/Variety/Abelian/Fq/2/5/a_g}{2.5.a\_g}.
    The particular overorder $S$ is the unique order with $[\cO_K:S]=8$, has $5$ weak equivalence classes and trivial Picard group.
    This gives $5$ isomorphism classes, with labels 
    $\texttt{2.5.a\_g-8.1.w.1}$ for $1\leq w\leq 5$.
    \item
    For $g=2$ and $q=5$, there are two isogeny classes (twists of each other) with maximum size of Picard group. 
    In both cases the Picard groups is $C_{12}$.
    See \href{https://abvar.lmfdb.xyz/Variety/Abelian/Fq/2/5/b_ac}{2.5.b\_ac} and 
    \href{https://abvar.lmfdb.xyz/Variety/Abelian/Fq/2/5/ab_ac}{2.5.ab\_ac}.
    \item
    For $g=2$ and $q=5$, whenever the Picard group has size $8$, it is isomorphic to $C_4 \times C_2$, and there are three such examples.  
    In one case, \href{https://abvar.lmfdb.xyz/Variety/Abelian/Fq/2/5/a_e}{2.5.a\_e}, the Frobenius order is maximal; the other two examples (\href{https://abvar.lmfdb.xyz/Variety/Abelian/Fq/2/5/b_e}{2.5.b\_e} and \href{https://abvar.lmfdb.xyz/Variety/Abelian/Fq/2/5/ab_e}{2.5.ab\_e}) have Frobenius order of index $50$.
\end{enumerate}

\newpage

\subsection*{Acknowledgements}
We thank the anonymous referees for corrections and suggestions.
Costa was supported by Simons Foundation grant SFI-MPS-Infrastructure-00008651.
Costa and Roe were supported by Simons Foundation grant 550033.
Dupuy is supported by National Science Foundation grant DMS-2401570.
Marseglia was supported by NWO grant VI.Veni.202.107,
and by Marie Sk{\l}odowska-Curie Actions - Postdoctoral Fellowships 2023 (project 101149209 - AbVarFq).
Vincent is supported by a Simons Foundation Travel Support for Mathematicians grant.

\printbibliography
\end{document}

\smm{rules to be followed: amsart, default margins, font size (10pt), and line spacing,  4 pages or less, excluding references. }